\documentclass[brazil]{article}
\usepackage{tgpagella}

\usepackage[T1]{fontenc}
\usepackage[latin9]{inputenc}
\usepackage{color}
\usepackage{babel}
\usepackage{longtable}
\usepackage{booktabs}
\usepackage{amsmath}
\usepackage{amssymb}
\usepackage{setspace}
\usepackage{subscript}
\usepackage[unicode=true,
 bookmarks=false,
 breaklinks=true,pdfborder={0 0 0},pdfborderstyle={},backref=page,colorlinks=true]
 {hyperref}
\hypersetup{
 pdfborderstyle=}

\makeatletter

\providecommand{\tabularnewline}{\\}

\@ifundefined{date}{}{\date{}}
\usepackage[alf]{abntex2cite}
\citeoption{abnt-etal-cite=2}
\usepackage{textcomp}
\usepackage{graphicx}
\usepackage{fancyhdr}
\providecommand{\tabularnewline}{\\}
\usepackage{syllogism}
\newcommand{\mh}{\noindent}
\hypersetup{citecolor=blue}

\makeatother

\begin{document}
\title{\textbf{\Large{}A definição de verdade de Tarski}\textbf{ }\thanks{Partes deste texto foram publicadas em Rodrigues, A. `Sobre a Concepção
da Verdade em Tarski', \textit{Abstracta}, v. 2, n. 1, p. 24-61, 2005.}\\\textbf{\Large{}Tarski's definition of truth}}
\author{{\large{}Guilherme Cardoso\thanks{Philosophy Department, Federal University of Ouro Preto (UFOP), guilherme.cardoso@ufop.edu.br.}   
 ~\& Abilio Rodrigues\thanks{Philosophy Department, Federal University of Minas Gerais (UFMG), abilio@ufmg.br, \url{https://orcid.org/0000-0001-6639-9550}.}}}
\maketitle
\begin{center}
\begin{minipage}[t]{0.8\columnwidth}%
\begin{center}
\textbf{\footnotesize{}Resumo}{\footnotesize\par}
\par\end{center}
{\footnotesize{}O objetivo deste texto é apresentar, de um modo tecnicamente
acessível, a definição de verdade de Tarski, o teorema da indefinibilidade
da verdade, e discutir duas críticas a aspectos conceituais do trabalho
de Tarski sobre a verdade, a saber, se a definição captura a noção
de verdade como correspondência e a objeção de Kripke à hierarquia
de linguagens.}\textbf{\footnotesize{}\\Palavras-chave:}{\footnotesize{}
Tarski, verdade, esquema T, indefinibilidade da verdade.}{\footnotesize\par}%
\end{minipage}
\par\end{center}

\begin{center}
\begin{minipage}[t]{0.8\columnwidth}%
\begin{center}
\textbf{\footnotesize{}Abstract}{\footnotesize\par}
\par\end{center}
{\footnotesize{}The aim of this text is to present, in a technically
accessible way, Tarski's definition of truth, the indefinability theorem,
and to discuss two aspects of Tarski's work on truth, namely, whether
or not the definition captures the notion of truth as correspondence,
and Kripke's objection to the hierarchy of languages.}\textbf{\footnotesize{}\\Keywords:}{\footnotesize{}
Tarski, truth, T-schema, indefinability of truth.}{\footnotesize\par}%
\end{minipage}
\par\end{center}

\section{Introdução}

Alfred Tarski foi um lógico e matemático polonês cuja {definição
de verdade}, apresentada nos anos 1930 na monografia \emph{O Conceito
de Verdade nas Linguagens Formalizadas} (\citeyear{Tarski2006a},
daqui em diante \emph{CVLF})\footnote{A monografia de Tarski sobre a verdade foi publicada pela primeira
vez em 1933 em polonês. Foi traduzida para o alemão em 1936 com o
acréscimo de um pós-escrito em que algumas teses de 1933 são revistas
e modificadas (ver \emph{CVLF,} p. 152, nota) e para o inglês (\citeauthor{Tarski1983}, \citeyear{Tarski1983}).
O texto base utilizado aqui é a tradução portuguesa que está na coletânea
\textit{A Concepção Semântica de Verdade} (\citeauthor{Tarski2006}, \citeyear{Tarski2006}).}, ocupa um lugar central nas discussões sobre a verdade, sobretudo
aquelas realizadas pela filosofia analítica. O objetivo de Tarski
em \emph{CVLF} é elaborar, para uma linguagem formalizada, uma definição
\emph{materialmente adequada} e \emph{formalmente correta} do predicado
verdade. A definição deveria também capturar a noção de verdade como
\emph{correspondência}, que Tarski chama de concepção clássica (\emph{CVLF,}
p. 19-20). O objetivo deste texto é apresentar, de um modo tecnicamente
acessível, 
a definição de verdade (Seção \ref{main.def}), o teorema da indefinibilidade
da verdade (Seção \ref{indef}), e discutir duas críticas a aspectos
conceituais do trabalho de Tarski, a saber, se a definição captura
a noção de verdade como correspondência (Seção \ref{corr}) e a objeção
de Kripke à hierarquia de linguagens (Seção \ref{kpk}). Ainda que
seja controverso se a definição expressa de modo satisfatório a noção
de verdade como correspondência, por outro lado, a concepção {deflacionista}
do problema da verdade, que se extrai do trabalho de Tarski, é uma
posição filosoficamente relevante acerca do problema da verdade.

\section{A definição de verdade}

\label{main.def}

\subsection{Adequação material e correção formal}

\label{amcf}

Uma definição de verdade {materialmente adequada} apresenta uma
condição que se aplica a todas e apenas às sentenças verdadeiras de
uma dada linguagem. Esse é um critério que diz respeito à {extensão}
da definição. Tarski propõe a \emph{Convenção T} como critério de
adequação material: uma definição de verdade em uma linguagem \emph{L}
será materialmente adequada se tiver como consequência todas as instâncias
do {esquema T}, 
\begin{description}
\item [{(T)}] a sentença \emph{S} é verdadeira se, e somente se, \emph{p}, 
\end{description}
\mh relativas às sentenças de \emph{L} (\emph{CVLF} p. 55-56). Instâncias
do esquema T são obtidas pela substituição de \emph{S} pelo nome de
uma sentença e de \emph{p} pela própria sentença, como por exemplo 
\begin{description}
\item [{(1)}] a sentença `Aristóteles é grego' é verdadeira se, e somente
se, Aristóteles é grego. 
\end{description}
\mh O esquema T não é uma definição de verdade, mas instâncias de
T são consideradas por Tarski definições parciais no sentido de definirem
verdade para as respectivas sentenças (\emph{CVLF} p. 23). (1) é uma
definição da verdade da sentença `Aristóteles é grego' porque expressa
a condição em que tal sentença é verdadeira, a saber, Aristóteles
ser grego.

Uma definição formalmente correta deve respeitar as regras para construção
de definições e as leis usuais da lógica\footnote{\emph{CVLF} p. 21 e 32-33. Ver também \emph{A concepção semântica
da verdade} (\citeauthor{Tarski2006c}, \citeyear{Tarski2006c}, daqui
em diante \emph{CSV}, p. 172): ``pressupomos que as regras formais
de definição são observadas na metalinguagem'' .}. Uma definição deve satisfazer a duas condições: (i) o \emph{definiendum}
deve ser eliminável de uma expressão em que ocorra, substituído pelo
\emph{definiens}; (ii) a definição nada pode acrescentar ao que já
`estava disponível' antes da introdução da definição -- isto é, a
definição, adicionada a uma teoria, não pode permitir a derivação
de um teorema que não fosse parte dessa teoria antes da introdução
da definição. Essas condições são denominadas, respectivamente, \emph{eliminabilidade}
e \emph{não criatividade}\footnote{Em \emph{CVLF} não há maiores esclarecimentos acerca das regras para
construção de definições, mas em uma nota de um texto de 1934 `Some
Methodological Investigations on the Definability of Concepts' (Tarski
1956 p. 307), Tarski menciona {não-criatividade} e {eliminabilidade}
como as duas condições que devem ser satisfeitas por uma definição
correta.}. Uma definição que implica uma contradição, sendo a lógica subjacente
clássica, claramente não satisfaz o critério de não criatividade,
pois a lógica clássica é explosiva, i.e. de uma contradição se segue
qualquer coisa. E as ``leis usuais da lógica'', para Tarski, são
as da lógica clássica.

Considerando que instâncias do esquema T são definições parciais de
verdade, para uma linguagem com um número {finito} de sentenças,
uma definição formalmente correta e materialmente adequada pode ser
obtida com base nas respectivas instâncias do esquema T. A título
de exemplo, suponha uma linguagem que possua apenas duas sentenças,
`a neve é branca' e `a grama é verde'. Uma definição de verdade para
essa linguagem seria dada pelo esquema abaixo: 
\begin{description}
\item [{(2)}] Para todo x, x é uma sentença verdadeira se, e somente se,
(x = `a neve é branca' e a neve é branca) ou (x = `a grama é verde'
e a grama é verde). 
\end{description}
\mh Ao considerar que instâncias do esquema T são definições parciais
de verdade, Tarski assume uma posição \emph{deflacionista} acerca
do problema da verdade. A ideia central do deflacionismo é que não
existe propriamente um problema filosófico acerca da natureza da verdade
e que tal noção é perfeitamente esclarecida na medida em que se reconhece
a equivalência entre a atribuição do predicado verdade a uma sentença
e a afirmação da própria sentença, precisamente o que é expressado
pelo {esquema T}.

\subsection{Contexto e motivações}

\label{ctx}

O trabalho de Tarski está de acordo com a ideia de dar à filosofia
um caráter científico, excluindo ingredientes considerados metafísicos.
Tarski pretendia reduzir conceitos semânticos a conceitos físicos,
lógicos e matemáticos. De fato, Tarski vai {eliminar} o predicado
verdade com o auxílio da equivalência expressada pelo esquema T. No
\emph{definiens} restarão apenas sentenças, conceitos lógicos e matemáticos,
e os conceitos aos quais a linguagem em questão se refere.

Em um texto de 1936, Tarski deixa claro que a semântica deve ser adequada
aos princípios da unidade da ciência e do fisicalismo \cite[p. 154]{Tarski2006b},
notoriamente defendidos pelo Círculo de Viena. Mas o mais provável
é que Tarski tenha sido influenciado não pelo Círculo de Viena, mas
sim pela Escola de Lvóv-Varsóvia, um importante movimento filosófico
polonês da primeira metade do século XX\footnote{Sobre os antecedentes e o contexto histórico do trabalho de Tarski,
ver \citeonline{Rojszczak2002}.}.

Essa visão que rejeita entidades metafísicas fica clara na opção de
Tarski por sentenças como portadores-de-verdade. Dos três possíveis
candidatos a portadores-de-verdade, crenças, proposições e sentenças,
apenas estas últimas podem ser analisadas em termos puramente físicos
e matemáticos\footnote{Ver \emph{CVLF} p. 23-24 (nota 3): ``é conveniente estipular que
termos como `palavra', `expressão', `sentença' etc. não denotam séries
concretas de sinais, mas a classe de todas aquelas séries cuja forma
é igual a da série dada''. Ver também \emph{CSV} p. 159.}. Segundo Tarski, sentenças não são inscrições particulares mas sim
conjuntos de inscrições com a mesma forma. Inscrições são objetos
físicos e conjuntos são objetos matemáticos.

A principal motivação de Tarski era mostrar que a noção de verdade
poderia ser usada de modo consistente em investigações lógicas. O
próprio Tarski observa que havia uma série de resultados que somente
poderiam ser adequadamente demonstrados se a noção de verdade fosse
precisamente definida (\emph{CVLF} p. 109-111, notas 83 e 85). Considere,
por exemplo, o teorema da completude da lógica clássica de primeira
ordem, 
segundo o qual toda {fórmula válida} pode ser demonstrada por meio
dos axiomas e regras de inferência do sistema dedutivo. Uma fórmula
válida é uma fórmula sempre verdadeira, ou verdadeira seja qual for
a interpretação e o domínio de discurso. Se a noção de verdade for
capaz de produzir contradições, esses resultados metateóricos que
lançam mão do conceito de verdade ficam sob suspeita.

\subsection{O paradoxo do mentiroso}

O chamado \emph{paradoxo do mentiroso} produz uma contradição a partir
de premissas aparentemente plausíveis e compromete o uso intuitivo
do predicado verdade. Dada uma sentença autorreferente como
\begin{description}
\item[] (S) a sentença S não é verdadeira, 
\end{description}
\mh supondo a validade irrestrita do esquema T e sendo clássica a
lógica subjacente, uma contradição é obtida nos seguintes passos: 
\begin{quote}
(3) a sentença `a sentença S não é verdadeira' é verdadeira se, e
somente se, a sentença S não é verdadeira.

(4) a sentença S é verdadeira se, e somente se, a sentença S não
é verdadeira.

(5) a sentença S é verdadeira ou a sentença S não é verdadeira.

Logo,

(6) a sentença S é verdadeira e a sentença S não é verdadeira. 
\end{quote}
No argumento acima, (3) é a instância do esquema T correspondente
à sentença S, (4) é obtida pela substituição da expressão `a
sentença S não é verdadeira' por S, nomes da mesma sentença,
e (5) é uma instância do princípio do terceiro excluído. Em ambos
os casos (tanto para S verdadeira quanto para S não verdadeira)
obtemos a contradição (6) (\emph{CVLF} p. 25)\footnote{A rigor, o Princípio do Terceiro Excluído não é necessário para derivar
uma contradição a partir de \emph{A} \ensuremath{\leftrightarrow}
\textasciitilde\emph{A}. Na lógica intuicionista, na qual não vale
o Terceiro Excluído, uma contradição também se segue de \emph{A} \ensuremath{\leftrightarrow}
\textasciitilde\emph{A.}}.

De acordo com o diagnóstico clássico, normalmente atribuído à Tarski,
são três os pressupostos que levam ao paradoxo\footnote{\emph{CVLF} p. 32 e \emph{CSV} p. 168-169.}: 
\begin{quote}
(I) A linguagem possui recursos para \emph{falar sobre} suas próprias
expressões e atribuir o predicado `\emph{x é verdadeira'} às suas
próprias sentenças. Tarski chama linguagens com essas características
de \emph{semanticamente fechadas}. Note que uma linguagem semanticamente
fechada é condição necessária para formular uma sentença autorreferente
como ({S}).

(II) As `leis usuais' da lógica valem -- ou seja, a lógica subjacente
é clássica, e portanto explosiva, tornando a teoria trivial na presença
de uma contradição.

(III) Todas as instâncias do esquema T são verdadeiras.\footnote{Na verdade, Tarski não apresenta explicitamente esses três pressupostos.
Tarski inclui o esquema T como parte do fechamento semântico (ver
\emph{CSV} p. 168). Entretanto, discussões posteriores tratam o problema
a partir da imposição de restrições à linguagem, à lógica ou ao esquema
T, precisamente os pressupostos acima (\cfcite{Feferman2008}).} 
\end{quote}
Dentre esses três pressupostos, pelo menos um deve ser rejeitado.
Tarski nem considera rejeitar (II) ou (III), rejeita (I) e conclui
que o paradoxo do mentiroso mostra que para uma linguagem semanticamente
fechada o predicado verdade produz uma contradição. Como a linguagem
natural é semanticamente fechada, dado que contém sua própria metalinguagem
e o predicado verdade aplicado às suas próprias expressões, Tarski
conclui que não é possível definir verdade para a linguagem natural
(\emph{CVLF} p. 33).

A solução proposta por Tarski será definir verdade para uma linguagem
\emph{L} que não tenha recursos para falar de sua própria semântica.
Tais linguagens são denominadas \emph{semanticamente abertas}. Assim,
verdade em \textit{L} será definida em uma metalinguagem \emph{ML}
de \emph{L} (CVLF p. 35). \emph{ML} deverá ter recursos para falar
da semântica de \emph{L} e atribuir o predicado verdade às sentenças
de \emph{L}. Mas \emph{ML} deverá também ser semanticamente aberta:
sua semântica somente poderá ser tratada na metametalinguagem, caso
contrário o paradoxo reaparece em \emph{ML}.

A solução proposta por Tarski, portanto, implica em uma \emph{hierarquia
de linguagens}, todas semanticamente abertas. O paradoxo do mentiroso
é evitado em uma linguagem semanticamente aberta, e que tenha seus
conceitos semânticos definidos em uma metalinguagem também semanticamente
aberta, porque o argumento usual que leva ao paradoxo não pode ser
reproduzido. Como o predicado verdade não pertence à linguagem objeto,
mas sim à metalinguagem (ou à linguagem de nível imediatamente superior),
não é possível formular uma sentença como ({S}), que diz de si mesma
que não é verdadeira.

\subsection{A definição de verdade para a linguagem do cálculo de classes}

Tarski elabora uma definição do predicado verdade para uma dada linguagem
\emph{L} que satisfaz certas condições. Uma dessas condições, se \emph{L}
tem um número infinito de sentenças, é que \emph{L} seja \emph{formalizada.}
Uma linguagem formalizada (também chamada de \emph{regimentada}) é
uma linguagem que tem a estrutura sintática precisamente estabelecida.
As sentenças de \emph{L} devem ter sido geradas {indutivamente}
a partir de expressões mais simples. Isso é condição necessária para
que o {valor semântico} de uma expressão complexa dependa unicamente
da sua estrutura e dos valores semânticos das partes que a compõem.

Um conjunto gerado indutivamente é um conjunto produzido a partir
da aplicação de um certo número de operações a um determinado conjunto
de elementos iniciais chamado \emph{base}.\footnote{O conjunto dos números naturais é gerado desse modo. A base é constituída
apenas do número 0, e todos os outros números são produzidos pela
operação \emph{sucessor}: 1 é o sucessor de 0, 2 é o sucessor de 1,
e assim por diante.} Utilizado na construção de uma linguagem, o método indutivo possibilita
produzir um número infinito de expressões a partir de um número finito
de operações aplicadas ao vocabulário também finito da linguagem.
O método indutivo permite também determinar se um predicado \emph{P}
(por exemplo, o predicado verdade) se aplica ou não aos elementos
de um conjunto \emph{C} gerado indutivamente (no caso, o conjunto
das expressões da linguagem), utilizando regras segundo as quais a
aplicação de \emph{P} a elementos de \emph{C} depende das condições
de aplicação de \emph{P} aos elementos mais simples de \emph{C} e
das operações utilizadas para gerar os demais elementos de \emph{C}.

Tarski não apresenta um método geral, mas sim uma definição de verdade
para uma determinada linguagem, formalizada e semanticamente aberta,
a \emph{linguagem do cálculo de classes} (\emph{LCC}). Note que Tarski
usa a palavra `linguagem' em um sentido mais forte que o usual, que
inclui um conjunto de axiomas e regras de inferência. A rigor, portanto,
uma linguagem é uma teoria. \emph{LCC} contém um sistema de lógica
de primeira ordem e axiomas do cálculo de classes.

Além das constantes lógicas, \emph{LCC} possui uma única constante
não lógica, o predicado binário \emph{I} que representa a relação
de inclusão. A fórmula \emph{Ix}\textsubscript{k}\emph{x}\textsubscript{l}
significa que \emph{x}\textsubscript{k} está contido em (ou é um
subconjunto de) \emph{x}\textsubscript{l}.
\begin{quote}
Vocabulário de \emph{LCC}:

1. operadores lógicos $\forall,{\sim},\lor$; 

2. o predicado binário \emph{I};

3. um número infinito de variáveis \emph{x}\textsubscript{k} (para
\emph{k} inteiro positivo);

4. sinais de pontuação `(' e `)'.\footnote{O vocabulário de \emph{LCC} apresentado acima é diferente do encontrado
em \emph{CVLF}. Este é um dentre alguns ajustes na terminologia e
nos símbolos utilizados que fazemos neste texto para torná-lo amigável
a um leitor que tenha familiaridade, por exemplo, com a lógica de
primeira ordem tal como é apresentada em \citeonline{Mortari2001}.} 
\end{quote}
Não há constantes individuais (i.e. nomes de indivíduos) em \emph{LCC}.
Os símbolos utilizados para as variáveis possibilitam a sua ordenação
em sequência, o que será essencial mais adiante.

A definição de verdade de \emph{LCC} será formulada em uma metalinguagem
de \emph{LCC}, que chamaremos de \emph{MLCC}. \emph{MLCC} deverá conter
nomes de todas as expressões de \emph{LCC} e expressões com o mesmo
significado das expressões de \emph{LCC}. Sendo \emph{e} uma expressão
da linguagem objeto \emph{LCC}, o nome de \emph{e} na metalinguagem
\emph{MLCC} será representado por \underline{\textit{e}}, e a expressão
da metalinguagem com o mesmo significado que \emph{e} será representada
por \underline{\underline{\textit{e}}}. Assim, a instância de T referente
à sentença
\begin{description}
\item [{(7)}] \textit{\ensuremath{\forall}x}\textsubscript{1}\textit{\ensuremath{\forall}x}\textsubscript{2}\textit{Ix}\textsubscript{1}\textit{x}\textsubscript{2} 
\end{description}
\mh de \emph{LCC} é
\begin{description}
\item [{(8)}] a sentença \textit{\underline{\ensuremath{\forall}x\textsubscript{1}\ensuremath{\forall}x\textsubscript{2}Ix\textsubscript{1}x\textsubscript{2}}}
é verdadeira se, e somente se, \textit{\underline{\underline{\ensuremath{\forall}x\textsubscript{1}\ensuremath{\forall}x\textsubscript{2}Ix\textsubscript{1}x\textsubscript{2}}}}. 
\end{description}
\mh Em (8), \underline{`{\ensuremath{\forall}x\textsubscript{1}\ensuremath{\forall}x\textsubscript{2}Ix\textsubscript{1}x\textsubscript{2}}'}
é um \emph{nome} de (7) e \underline{\underline{`\ensuremath{\forall}x\textsubscript{1}\ensuremath{\forall}x\textsubscript{2}Ix\textsubscript{1}x\textsubscript{2}'}}
é a \emph{tradução} de (7) na metalinguagem \emph{MLCC}.

\medskip{}

\begin{quote}
Definição de fórmula de \emph{LCC}:

1. Para \emph{k} e \emph{l} inteiros positivos, \underline{$Ix_l x_k$}
é fórmula;

2. Se \emph{A} é fórmula, \underline{${\sim}$} $A$ é fórmula;

3. Se \emph{A} e \emph{B} são fórmulas, \emph{A}\underline{{\ensuremath{\vee}}}\emph{B}
é fórmula;

4. se \emph{A} é fórmula, \underline{{\textit{\ensuremath{\forall}x\textsubscript{k}}}}\emph{A}
é fórmula, para \emph{k} inteiro positivo;

5. Nada mais é fórmula.\footnote{Ver Definição 10 \emph{CVLF} p. 44. Note que uma definição recursiva
como essa não atende o critério de eliminabilidade. Tarski na nota
24 (p. 44-45) apresenta a definição normal equivalente, na qual o
\emph{definiendum} não ocorre no \emph{definiens.}} 
\end{quote}
As letras \emph{A} e \emph{B} são variáveis da metalinguagem que varrem
fórmulas de \emph{LCC}. Uma sentença, como usual, é uma fórmula sem
variáveis livres, i.e. uma fórmula fechada.

Note que todas as fórmulas atômicas de \emph{LCC} são abertas (cláusula
1). Por essa razão, a verdade de sentenças complexas não pode ser
definida a partir da verdade de sentenças mais simples.\footnote{Para uma linguagem sem fórmulas abertas, por exemplo, uma linguagem
da lógica sentencial, é possível definir verdade de sentenças complexas
a partir da verdade das sentenças mais simples\emph{.} Lembre como
as tabelas de verdade estabelecem o valor de verdade de uma sentença
complexa a partir dos valores de verdade das sentenças atômicas.} A sentença (7) tem as seguintes subfórmulas: (i) \textit{Ix}\textsubscript{1}\textit{x}\textsubscript{2};
(ii) \textit{\ensuremath{\forall}x}\textsubscript{2}\textit{Ix}\textsubscript{1}\textit{x}\textsubscript{2}
e (iii) \textit{\ensuremath{\forall}x}\textsubscript{1}\textit{\ensuremath{\forall}x}\textsubscript{2}\textit{Ix}\textsubscript{1}\textit{x}\textsubscript{2}.
(i) e (ii) são fórmulas abertas, que não são nem verdadeiras nem falsas,
e (iii) é a própria sentença (7). Precisamos de uma noção mais geral
que verdade, que se aplique também a fórmulas abertas. Esse é o papel
da relação de \emph{satisfação} entre fórmulas da linguagem e objetos
do domínio -- ou mais precisamente, como veremos mais adiante, sequências
infinitas de objetos. O esquema 
\begin{description}
\item [{(9)}] \emph{a} satisfaz a fórmula \emph{F} se, e somente se, \emph{p} 
\end{description}
\mh funciona de maneira análoga ao esquema T. \emph{a} e \emph{F}
são substituídos respectivamente por nomes da metalinguagem de um
objeto e de uma fórmula, e \emph{p} é substituído por uma sentença
da metalinguagem que expressa a condição para que o objeto \emph{a}
satisfaça a fórmula \emph{F}. A título de exemplo, supondo que Aristóteles
e Descartes sejam elementos do domínio, (9) produz as seguintes sentenças
da metalinguagem relativas à fórmula aberta (ou predicado) `\emph{x}\textsubscript{1}\emph{
é grego}': 
\begin{description}
\item [{(9a)}] Aristóteles satisfaz `\emph{x}\textsubscript{1} é grego'
se, e somente se, \emph{Aristóteles é grego}; 
\item [{(9b)}] Descartes satisfaz `\emph{x}\textsubscript{1} é grego'
se, e somente se, \emph{Descartes é grego}. 
\end{description}
\mh Uma fórmula aberta pode ter mais de uma variável livre. O predicado
`\emph{x}\textsubscript{1} é grego' é satisfeito por Aristóteles,
mas não por Descartes; o predicado `\emph{x}\textsubscript{1} é discípulo
de \emph{x}\textsubscript{2}' é satisfeito pelo par ordenado {[}Platão,
Sócrates{]}, mas não pelo par ordenado {[}Platão, Aristóteles{]}.

Para aplicar a noção de satisfação a \emph{LCC} de modo uniforme para
fórmulas com um número arbitrário de variáveis, Tarski usa sequências
infinitas de objetos (ver \emph{CVLF} p. 59). Uma sequência infinita
de objetos é uma atribuição de valores às variáveis da linguagem.
Do ponto de vista matemático, é uma função que vai do conjunto de
variáveis, indexadas por um inteiro positivo \emph{k}, ao universo
de discurso da linguagem. As sequências são infinitas porque a cada
variável indexada por um número inteiro positivo é atribuído um elemento
do domínio. Note que a cada posição \emph{k} vai corresponder um único
elemento do domínio, mas um mesmo elemento do domínio pode ocupar
mais de uma posição \emph{k}. A título de exemplo, suponha que o domínio
seja o conjunto

\smallskip{}

D = \{Sócrates, Platão, Aristóteles, Kant\}.

\smallskip{}

\mh As variáveis são indexadas por números inteiros positivos. As
quatro sequências infinitas representadas na tabela abaixo exemplificam
como podemos dispor as primeiras cinco posições (variáveis \emph{x}\textsubscript{1}
a \emph{x}\textsubscript{5}):\newpage

\begin{longtable}[c]{@{}lllllll@{}}
\toprule 
 & \emph{x}\textsubscript{1}\emph{ }  & \emph{x}\textsubscript{2}  & \emph{x}\textsubscript{3}\emph{ }  & \emph{x}\textsubscript{4}\emph{ }  & \emph{x}\textsubscript{5}\emph{ }  & ...\tabularnewline
\emph{f}\textsuperscript{1}  & Sócrates  & Platão  & Aristóteles  & Platão  & Platão  & ...\tabularnewline
\emph{f}\textsuperscript{2}  & Sócrates  & Sócrates  & Kant  & Aristóteles  & Kant  & ...\tabularnewline
\emph{f}\textsuperscript{3}  & Aristóteles  & Platão  & Sócrates  & Kant  & Sócrates  & ...\tabularnewline
\emph{f}\textsuperscript{4}  & Platão  & Aristóteles  & Platão  & Aristóteles  & Aristóteles  & ...\tabularnewline
\bottomrule
\end{longtable}

\mh Qualquer que seja o número de variáveis livres de uma fórmula
\emph{F}, a relação de satisfação se dá entre \emph{F} e sequências
infinitas, mas somente são considerados os objetos que ocupam posições
correspondentes às variáveis livres de \emph{F}, os outros são desprezados.
Uma sequência \emph{f} satisfaz a fórmula aberta \emph{Rx}\textsubscript{i}\emph{x}\textsubscript{j}
se o objeto da posição \emph{i} na sequência \emph{f} está na relação
\emph{R} com o objeto da posição \emph{j} da sequência \emph{f.} Agora
considere as fórmulas \emph{F}\textsubscript{1} = `\emph{x}\textsubscript{1}
é mestre de \emph{x}\textsubscript{2}', e \emph{F}\textsubscript{2}
= `\emph{x}\textsubscript{4} não é grego'. As sequências \emph{f}\textsuperscript{1}
e \emph{f}\textsuperscript{4} satisfazem \emph{F}\textsubscript{1},
mas não satisfazem \emph{F}\textsubscript{2}. \emph{F}\textsubscript{2}
é satisfeita apenas pela sequência de objetos \emph{f}\textsubscript{3}.

Em \emph{MLCC} (a metalinguagem na qual verdade será definida) o \emph{k}-ésimo
elemento de uma sequência \emph{f}, denotado por \emph{f}\textsubscript{k},
é o valor atribuído à variável \emph{x}\textsubscript{k}. Tarski
formula então o seguinte esquema (\emph{CVLF} p. 60):
\begin{description}
\item [{(10)}] a sequência infinita \emph{f} satisfaz a fórmula \emph{F}
se, e somente se, \emph{p}. 
\end{description}
\mh A instância de (10) relativa a uma fórmula \emph{F} e uma dada
sequência \emph{f} é obtida substituindo-se \emph{F} por um nome de
\emph{F} na metalinguagem, e no lugar de \emph{p} coloca-se a expressão
obtida pela substituição das variáveis livres de \emph{F} pelos símbolos
`\emph{f}\textsubscript{k}', `\emph{f}\textsubscript{l}', etc. que
denotam os objetos que na sequência \emph{f} correspondem às posições
\emph{k}, \emph{l}, etc. das variáveis livres de \emph{F.} A instância
de (10) relativa à fórmula \textit{Ix}\textsubscript{1}\textit{x}\textsubscript{2}
é 
\begin{description}
\item [{(11)}] \emph{f} satisfaz a fórmula \textit{\underline{Ix\textsubscript{1}x\textsubscript{2}}}
se, e somente se \underline{\underline{\textit{I}}}\emph{f}\textsubscript{1}\emph{f}\textsubscript{2}. 
\end{description}
\mh Considere uma sequência \emph{f'} com \ensuremath{\varnothing}
e \{\{\ensuremath{\varnothing}\}\} respectivamente nas posições 1
e 2. `\emph{f}\textsubscript{1}' e `\emph{f}\textsubscript{2}' serão
nomes da metalinguagem respectivamente para \ensuremath{\varnothing}
e \{\{\ensuremath{\varnothing}\}\}. \emph{f'} satisfaz a fórmula
$Ix_{1}x_{2}$ (lembre-se que o conjunto vazio é subconjunto de todos
os conjuntos). Por outro lado, uma sequência \emph{f''} com \{\ensuremath{\varnothing}\}
e \{\{\ensuremath{\varnothing}\}\} respectivamente nas posições 1
e 2 não satisfaz a fórmula $Ix_{1}x_{2}$.

A linguagem de \emph{LCC} foi definida indutivamente. Agora, vamos
definir recursivamente a relação de satisfação. As condições segundo
as quais uma sequência de objetos \emph{f} satisfaz uma fórmula atômica
são estabelecidas diretamente pela cláusula de base e, em seguida,
é estabelecido como \emph{f} se comporta em relação às operações utilizadas
na construção das fórmulas complexas de \emph{LCC}.

\medskip{}

\begin{quote}
Definição de satisfação:

Seja \emph{f} um sequência de objetos e $f_{i}$ 
(para \emph{i} inteiro positivo) o nome do \emph{i}-ésimo elemento
de \emph{f}.

1. \emph{f} satisfaz \textit{{Ix}\textsubscript{k}\textit{x}\textsubscript{l}\textit{}}
sse \underline{\underline{\it I}}\textit{f}\textsubscript{k}\textit{f}\textsubscript{l},
para k e l inteiros positivos;

2. \emph{f} satisfaz \underline{\textasciitilde}\emph{A} sse \emph{f}
não satisfaz \emph{A};

3. \emph{f} satisfaz \emph{A} {\underline{\ensuremath{\vee}}} \emph{B}
sse \emph{f} satisfaz \emph{A} ou \emph{f} satisfaz \emph{B};

4. \emph{f} satisfaz \underline{\it{\ensuremath{\forall}x\textsubscript{k}}}\emph{A}
se e somente se, toda sequência que difere de \emph{f} no máximo em
seu k-ésimo lugar satisfaz \emph{A}.\footnote{\emph{CVLF} p. 61, Definição 22. Na nota 41, a definição normal equivalente.} 
\end{quote}
\mh Note que cada cláusula da definição de satisfação corresponde
a uma cláusula da definição de fórmula. A linguagem ter uma sintaxe
precisamente determinada (ou ter sido gerada indutivamente, ou ser
\emph{formalizada}) é uma condição necessária para que a noção de
satisfação possa ser também precisamente definida.

Note que tanto a definição de satisfação quanto a de fórmula não atendem
o critério de eliminabilidade. Tarski transforma essas definições
em definições explícitas (que ele chama de \emph{normais}) usando
recursos da teoria de conjuntos, de modo a eliminar do \emph{definiens}
todas as ocorrências do \emph{definiendum}. A definição explícita
substitui os termos semânticos que ocorrem no \emph{definiens} por
termos lógicos e matemáticos, o que atende, portanto, a exigência
de reduzir a semântica a noções físicas, lógicas e matemáticas. As
definições explícitas (ou normais) de fórmula e satisfação são apresentadas
nas notas 24 e 41 de \emph{CVLF}.

No caso de sentenças, há apenas duas possibilidades: ou a sentença
é satisfeita por todas as sequências ou por nenhuma. No primeiro caso
a sentença é verdadeira e no segundo é falsa. Chegamos à definição
de verdade de \emph{LCC} (\emph{CVLF} p. 63):
\begin{description}
\item [{(V)}] para toda sentença \emph{x}, \emph{x} é verdadeira se, e
somente se, para toda sequência \emph{f}, \emph{f} satisfaz \emph{x}. 
\end{description}
\mh Vamos ilustrar como funciona a cláusula 4 da definição de satisfação
com as sentenças
\begin{description}
\item [{(12)}] Para todo \emph{x}\textsubscript{1}, \emph{x}\textsubscript{1}
é filósofo, 
\item [{(13)}] Para todo \emph{x}\textsubscript{1}, \emph{x}\textsubscript{1}
é grego, 
\end{description}
\mh e o domínio
\begin{description}
\item [{D}] = \{Sócrates, Platão, Aristóteles, Kant\}.
\end{description}
\mh Claramente, as sentenças (12) e (13) são, respectivamente, verdadeira
e falsa em relação a esse domínio de indivíduos. Segundo a cláusula
4 da definição de satisfação, uma dada sequência \emph{f} satisfaz
a sentença (12) se e somente se toda sequência que difere de \emph{f}
no máximo na posição 1 satisfaz a fórmula aberta `\emph{x}\textsubscript{1}
é filósofo'. Note que qualquer sequência de objetos, inclusive \emph{f},
terá na posição 1 um indivíduo que é filosófo (só há filósofos no
domínio). Portanto, todas as sequências satisfazem a sentença (12).
No caso da sentença (13) isso não ocorre, pois há sequências com Kant
na posição 1, e Kant não é grego. Portanto, nenhuma sequência satisfaz
(13), que por essa razão é falsa.

Na linguagem utilizada por Tarski em \emph{CVLF} não há constantes,
i.e. nomes de indivíduos. Mas cabe perguntar como a definição de verdade
em termos de satisfação funcionaria para uma linguagem de primeira
ordem com constantes. Considere a sentença 
\begin{description}
\item [{(14)}] Aristóteles é grego. 
\end{description}
\mh Não há variáveis livres e nem quantificadores em (14). Nesse
caso, a condição para que uma sequência \emph{f} satisfaça a sentença
é dada pela própria sentença: 
\begin{description}
\item [{(14')}] a sequência \emph{f} satisfaz `Aristóteles é grego' se,
e somente se, Aristóteles é grego. 
\end{description}
\mh Como Aristóteles é de fato grego, a condição à direita da bicondicional
é satisfeita por todas as sequências, de acordo com a definição de
verdade em termos de satisfação. E de fato, (14) é verdadeira.

\section{O Teorema da Indefinibilidade da Verdade}

\label{indef}

O método ilustrado acima para a construção de definições materialmente
adequadas e formalmente corretas dos predicados de verdade de certas
linguagens é normalmente reconhecido como a parte positiva dos resultados
estabelecidos por Tarski em \emph{CVFL}. Existem, entretanto, linguagens
(ou melhor, teorias) para as quais o método de Tarski não funciona.
Isso, grosso modo, é o que nos diz o Teorema da Indefinibilidade da
Verdade, apresentado em \emph{CVLF} (p. 117). O Teorema da Indefinibilidade
é reconhecido como a parte negativa dos resultados de Tarski, um teorema
de limitação.

Em \emph{CVLF}, Tarski enuncia o Teorema da Indefinibilidade para
o caso específico da \emph{Linguagem do Cálculo Geral de Classes}
(daqui em diante \emph{LCGC}). \emph{LCGC} pode ser obtida a partir
do vocabulário da \emph{LCC}, substituindo-se o esquema de variáveis
$x_{k}$ (onde \emph{k} é qualquer inteiro positivo) pelo esquema
de variáveis $x_{k}^{m}$ (onde \emph{m} e \emph{k} são inteiros positivos).
O índice subscrito \emph{k} permite a ordenação das variáveis, enquanto
o índice sobrescrito \emph{m} expressa a \emph{ordem} de cada variável.

Indivíduos são de ordem 1, classes de indivíduos são de ordem 2, classes
de classes de indivíduos são de ordem 3, assim por diante. Enquanto
a \emph{LCC} é uma linguagem de ordem finita (de ordem 2, pois suas
variáveis varrem classes, indistintamente), a \emph{LCGC} é uma linguagem
de ordem infinita.

O Teorema da Indefinibilidade diz basicamente, neste caso, que a classe
das sentenças verdadeiras da \emph{LCGC} não pode ser corretamente
definida na \emph{MLCGC} (a metalinguagem de \emph{LCGC}), de modo
que todas as sentenças da \emph{MLCGC} satisfaçam o esquema T e a
lógica subjacente seja clássica. O que ocorre aqui é que \emph{LCGC}
é `suficientemente rica' para expressar sua própria metalinguagem
(\emph{MLCGC}), produzindo instâncias do esquema T para todas as sentenças
de \emph{MLCGC}. Ao mesmo tempo, em \emph{LCGC} é possível construir
sentenças autorreferentes como {(S)}. Na presença da lógica clássica,
como já vimos, daí se segue uma contradição.

Uma versão mais conhecida do resultado acima apresentado, o Teorema
da Indefinibilidade da Verdade Aritmética, diz que \emph{a verdade
aritmética não pode ser definida na própria aritmética}. Apresentaremos
a seguir um esboço deste resultado.\footnote{O leitor interessado encontrará esse resultado apresentado em detalhe
em \citeonline[cap. 17]{Boolos2012}. Outra apresentação pode ser
encontrada em \citeonline{Cardoso2018}.}

Em um sistema formal da aritmética (i.e. axiomas da aritmética mais
um sistema de lógica de predicados de primeira ordem -- e.g. \citeonline[cap. 17]{Mortari2001})
é possível `falar sobre' sentenças desse próprio sistema formal. Vamos
chamar de \emph{Arit}\footnote{As teorias aritméticas que nos interessam são as extensões da Aritmética
de Robinson ($Q$). Nestas teorias se pode representar todas as funções
recursivas primitivas, portanto, podemos aplicar a elas o Lema Diagonal.} esse sistema formal e \emph{L}\textsubscript{A} a respectiva linguagem.
Pelo resultado denominado \emph{Lema da Diagonal}, para qualquer predicado
$P(x)$ da linguagem $L_{A}$, existe uma sentença $G$ de $L_{A}$
que `diz de si mesma' que ela satisfaz $P(x)$. Mais precisamente:
\begin{description}
\item [{$Arit\vdash G\leftrightarrow P(\#G)$}]~
\end{description}
\mh Onde $\#G$ é o numeral do número de gödel de $G$, ou seja,
um termo de $L_{A}$ que nomeia o número que codifica a sentença $G$.
Esse nome é obtido pelo método denominado \emph{aritmetização da sintaxe},
ou \emph{numeração de} Gödel, que atribui um número inteiro positivo
a cada expressão de $L_{A}$, um único número para cada expressão.
Em seguida, vamos ilustrar o funcionamento dessa codificação.

Em primeiro lugar, atribuímos um inteiro positivo para cada símbolo
primitivo de $L_{A}$. Segundo um importante resultado da aritmética,
que certamente o leitor vai recordar, todo inteiro positivo tem uma
única decomposição em fatores primos. Assim, podemos codificar uma
expressão qualquer de $L_{A}$ (a sequência de símbolos a ela associada)
pelo número cuja decomposição em fatores primos revela a mesma sequência
na posição dos exponenciais. Por exemplo, suponha que codificamos
os símbolos

\[
(,),\sim,0,=,s
\]

\noindent respectivamente, pelos inteiros

\[
1,3,7,8,13,24
\]

\mh Assim, a expressão

\[
\sim(0=s(0))
\]

\noindent será codificada pelo número cuja decomposição em fatores
primos apresenta tal sequência nas posições exponenciais. Reservamos
o expoente do número $2$ para informar o comprimento da sequência
(o número de caracteres, que na expressao `$\sim(0=\emph{s}(0))$'
é igual a 9). Assim, podemos codificar a expressão anterior por

\[
2^{9}\cdot3^{7}\cdot5^{1}\cdot7^{8}\cdot11^{13}\cdot13^{24}\cdot17^{1}\cdot19^{8}\cdot23^{3}\cdot29^{3}
\]

\mh Temos, portanto, uma maneira de codificar fórmulas de $L_{A}$
por meio de números. O nome de uma expressão em \emph{L}\textsubscript{A}
será o numeral do número de Godel que codifica tal expressão. Como
o leitor certamente notou, à toda expressão de $L_{A}$ irá corresponder
um único número de Gödel. Por outro lado, há inteiros positivos que
não são números de expressões `bem-formadas'\footnote{Por exemplo, ao número $2^{3}\cdot3^{3}\cdot5^{13}\cdot7^{1}$ corresponde
a expressão `$)=($', que claramente é apenas uma sequência de símbolos,
mas não uma expressão `bem-formada' da linguagem da aritmética.}. Mas dado um inteiro positivo qualquer, é possível determinar se
ele é o número de alguma expressão e, caso o seja, determinar qual
é a expressão correspondente de $L_{A}$.

Agora suponha que acrescentemos à linguagem $L_{A}$ um predicado
verdade, $V(x)$. Vamos chamar essa nova linguagem de $L_{A}^{*}$.
Seja $Arit^{*}$ uma teoria aritmética na linguagem $L_{A}^{*}$\footnote{Uma teoria que seja ao menos tão forte quanto $Q$, como dissemos
antes.} capaz de provar todas as instâncias do esquema T de $L_{A}^{*}$
-\/- i.e. todas as instâncias do esquema T valem para a linguagem
$L_{A}^{*}$. O predicado ${\sim}V(x)$ claramente pertence à linguagem
$L_{A}$. Portanto, pelo Lema Diagonal, existe uma sentença $M$ de
$L_{A}^{*}$, tal que:

\[
Arit^{*}\vdash M\leftrightarrow\ {\sim}V(\#M)
\]

\mh Mas como $L_{A}$ prova todas as instâncias do esquema T,

\[
Arit^{*}\vdash M\leftrightarrow V(\#M)
\]

\mh Uma contradição se segue em poucos passos.

Trocando em miúdos, uma teoria \textit{T} que tenha o poder expressivo
da aritmética dos números naturais (e é bastante razoável que qualquer
teoria minimamente interessante tenha o poder expressivo da aritmética
dos números naturais) e cuja linguagem possua o predicado verdade
será inconsistente (i.e., tal teoria produz uma contradição).

\section{O esquema T e a noção de verdade como correspondência}

\label{corr}

Embora a importância do trabalho de Tarski seja amplamente reconhecida,
por outro lado, há reações negativas. Nesta seção, abordaremos a seguinte
questão: a definição de Tarski expressa de modo satisfatório a noção
de verdade como correspondência?

A intuição básica da noção de verdade como correspondência é que uma
sentença é verdadeira em virtude de algo na realidade que funciona
como seu `fazedor-de-verdade' (\emph{truthmaker}). Tarski pretendia
que sua definição captasse a noção de verdade como correspondência,
pois isso é dito textualmente em \emph{CVLF}\footnote{\emph{CVLF} p. 153. Ver também \emph{SCT} p. 336 (\S\ 5).}.
Em \emph{CSV} (1944, p. 160) Tarski menciona novamente o trecho da
\textit{Metafísica} de Aristóteles, 
\begin{quote}
Dizer do que é que não é, ou do que não é que é, é falso, enquanto
que dizer do que é que é, ou do que não é que não é, é verdadeiro
(\textit{Metafísica} 1011b26-28). 
\end{quote}
como uma expressão da noção de verdade como correspondência (CSV,
p. 160). Note como o trecho acima é bastante próximo tanto da ideia
expressada pelo esquema T como também da relação de satisfação. A
analogia entre o trecho acima e a relação de satisfação é ainda mais
clara quando lembramos que, para Aristóteles, toda proposição é um
`dizer algo sobre algo' (\textit{Da Interpretação} 17a25) o que, no
caso de uma linguagem de primeira ordem, é atribuir um predicado a
uma \emph{n-}upla ordenada de objetos.

\citeonline[p. 71]{Putnam1979} 
dirige uma forte crítica ao trabalho de Tarski.
Segundo Putnam, instâncias de T dizem apenas que aceitar que \emph{p}
é verdadeira implica em aceitar \emph{p}, e vice-versa. O problema
é que isso seria muito pouco para adeptos da noção de verdade como
correspondência. Putnam menciona a última frase do romance nonsense
de Lewis Carrol, \emph{The Hunting of the Snark}: `the Snark was a
Boojum', que pretende significar nada além de expressar que o Snark
(personagem do romance) era um Boojum -- seja lá o que isso signifique.
De fato, é duvidoso que algum interessado na noção de verdade como
correspondência esteja interessado em tais instâncias do esquema T.

\citeonline[p. 314, 323]{Popper1972}, por outro lado, defende vigorosamente
o trabalho de Tarski sobre a verdade. Segundo Popper, Tarski reabilitou
a noção de verdade como correspondência ao apontar que a verdade de
uma linguagem \emph{L} deve ser tratada em uma metalinguagem de \emph{L},
pois para que possamos falar da relação entre sentenças e fatos (ou
o que quer seja que torne sentenças verdadeiras), precisamos falar
simultaneamente sobre a linguagem (objeto) e a realidade, e isso deve
ser feito na metalinguagem.

Popper está correto ao afirmar que foi um mérito de Tarski ter mostrado
que a relação entre linguagem (objeto) e mundo, central para a noção
de verdade como correspondência, deve ser expressada em uma metalinguagem
na qual se possa falar de ambos (linguagem objeto e mundo). Popper,
entretanto, não diz explicitamente por que instâncias do esquema T
expressam uma relação entre a linguagem e a realidade. Essa justificativa
é simples, embora para alguns possa não parecer convincente. Gila
Sher a apresenta da seguinte forma:
\begin{quote}
O ponto central do esquema T, do ponto de vista da correspondência,
é o contraste entre os lados esquerdo e direito de suas instâncias
(...) O lado esquerdo de uma bicondicional T é uma predicação linguística,
já o lado direito é uma predicação \emph{objetual}, que `diz respeito
ao mundo'. A tarefa de uma teoria da verdade como correspondência
é reduzir predicações de \emph{verdade}, que são linguísticas, a predicações
\emph{objetuais}. \footnote{\citeauthor{Sher1999b}, \citeyear{Sher1999b}, p. 135-6.}
\end{quote}
\mh Assim, a sentença da metalinguagem
\begin{description}
\item [{(1)}] a sentença `Aristóteles é grego' é verdadeira se, e somente
se, Aristóteles é grego, 
\end{description}
\noindent reduz uma predicação linguística, a atribuição do predicado
verdade à sentença `Aristóteles é grego', a uma predicação objetual,
a atribuição do predicado `\emph{x} é grego' a Aristóteles. O lado
direito diz respeito a um fato do mundo, contingente. Isso é possível
porque a metalinguagem tem recursos para falar das expressões da linguagem
objeto, atribuindo a tais expressões o predicado verdade. Além disso,
dado que tem expressões com o mesmo significado das expressões da
linguagem objeto, a metalinguagem tem também recursos para `falar
do mundo', ou seja, falar das mesmas entidades não linguísticas acerca
das quais a linguagem objeto fala. Em outras palavras, na linguagem
objeto há uma relação linguagem-mundo, ao passo que na metalinguagem
há tanto uma relação linguagem-mundo quanto uma relação linguagem-linguagem.
Por esse motivo, instâncias do esquema T não expressam apenas uma
relação entre uma sentença e seu nome, mas sim entre a sentença mencionada
no lado esquerdo e o mundo.

Outro argumento que defende a tese de que a definição de Tarski é
uma teoria de correspondência baseia-se na relação de satisfação.
A relação de satisfação é formulada por Tarski em \emph{CVLF}, por
motivos técnicos, como uma relação entre sequências infinitas de objetos
e fórmulas, mas é essencialmente uma relação entre objetos do universo
de discurso e fórmulas abertas. Por esse motivo, é interpretada como
um tipo de relação de correspondência na qual os objetos têm a função
de `\emph{truthmakers}' (\cfcite{Kolar1999}).

Vimos que satisfação é definida recursivamente sobre a complexidade
das fórmulas da linguagem para a qual está sendo definido o predicado
verdade. A relação propriamente dita entre objetos e fórmulas é dada
pela cláusula de base, que claramente funciona de modo análogo ao
esquema T. O papel da relação de satisfação é tornar possível a generalização
da idéia básica do esquema T para linguagens com um número infinito
de sentenças ou que não tenham nomes para todos os indivíduos do universo
de discurso. A razão pela qual Tarski define verdade usando satisfação
é puramente técnica. A rigor, não há uma prevalência da noção de satisfação
em relação à noção de verdade, nem vice-versa. Já vimos acima que
\begin{description}
\item [{(9a)}] Aristóteles satisfaz `\emph{x}\textsubscript{1} é grego'
se, e somente se, \emph{Aristóteles é grego}; 
\item [{(9b)}] Descartes satisfaz `\emph{x}\textsubscript{1} é grego'
se, e somente se, \emph{Descartes é grego}. 
\end{description}
\mh Assim como instâncias do esquema T, (9a) e (9b) não dizem se
Aristóteles e Descartes satisfazem ou não a fórmula `\emph{x}\textsubscript{1}
é grego', mas dizem apenas que sempre que se aceita ou rejeita o lado
esquerdo da bicondicional, deve-se também aceitar ou rejeitar o lado
direito.

Satisfação é uma relação linguagem mundo tanto quanto uma instância
do esquema T. As mesmas acusações que são dirigidas ao esquema T podem
ser também dirigidas à definição de satisfação. Por outro lado, (9a)
e (9b), de modo análogo a instâncias de T, podem ser interpretadas
como uma redução de uma predicação linguística a uma predicação objetual.

Mas será que os argumentos acima expostos são suficientes para que
a definição de Tarski seja bem sucedida na tentativa de captar a noção
de verdade como correspondência? Essa é uma questão que não vamos
responder aqui, mas sim deixar que o leitor tire suas próprias conclusões,
ou se aprofunde no assunto para tentar preencher lacunas, nos argumentos
a favor ou contra Tarski.

\section{A objeção de Kripke à hierarquia de linguagens}

\label{kpk}

Como vimos anteriormente, o diagnóstico clássico identifica o paradoxo
como resultado de três pressupostos: o fechamento semântico da linguagem
em questão, a lógica clássica e o esquema T. Os dois últimos pressupostos,
entretanto, não são negociáveis para Tarski\footnote{Ver \emph{CSV,} p. 169.}.
Logo, o paradoxo resulta do fechamento semântico. Esse resultado,
todavia, contém lições profundamente distintas a respeito da linguagem
natural e das linguagens formalizadas.

Em relação às linguagens formalizadas, o diagnóstico clássico pode
ser interpretado de uma perspectiva normativa, na medida em que estabelece
limites sobre o que podemos consistentemente expressar em tais linguagens.
Isso ocasiona a construção de hierarquias de linguagens que obedeçam
a tais limites de expressividade, impedindo o fechamento semântico.
Assim, os paradoxos são evitados de maneira eficiente mas, segundo
algumas análises, extremamente artificial e técnica, já que aparentemente
não há justificativa para tais restrições senão apenas evitar os paradoxos.

Quanto à linguagem natural, o diagnóstico não pode assumir uma perspectiva
normativa. A linguagem natural é semanticamente fechada e nada há
que se possa restringir (ver CVLF, p. 32). De fato, parece difícil
conceber exatamente como o diagnóstico e as hierarquias poderiam ser
implementados para que se pudesse extrair qualquer tipo de lição a
respeito daquilo que os paradoxos realmente revelam sobre o conceito
de verdade. Esta é a principal crítica levantada ao diagnóstico clássico.

Além disso, pode-se levantar uma objeção quanto ao caráter \emph{ad
hoc}, excessivamente intrínseco e artificial das hierarquias de linguagem.
Esta objeção é melhor explicada por meio de uma versão do paradoxo
formulada por Kripke\footnote{\cfcite[p. 54-55]{Kripke1984}.}. Considere
as seguintes sentenças:
\begin{description}
\item [{(15)}] A maioria das asserções de Nixon sobre Watergate são falsas. 
\item [{(16)}] Tudo aquilo que Jones diz sobre Watergate é verdadeiro. 
\end{description}
Suponha adicionalmente que (15) é a única sentença asserida por Jones
a respeito do caso Watergate, que Nixon asseriu (16) e que todas as
sentenças asseridas por Nixon acerca do Watergate, excluindo-se apenas
(16), são igualmente divididas entre aquelas que são verdadeiras e
aquelas que são falsas. Assim, o leitor pode facilmente verificar
que (15) é verdadeira sse (15) é falsa e, pelas mesmas razões, (16)
é verdadeira sse (16) é falsa.

Em primeiro lugar, a hierarquia de linguagens não diz nada sobre o
paradoxo de Kripke, já que (15) e (16) são sentenças de uma mesma
linguagem. Algumas soluções de inspiração tarskiana, entretanto, tem
defendido que a hierarquia de linguagens de Tarski poderia ser substituída
por uma hierarquia de predicados semânticos e que isto permitiria
solucionar os paradoxos\footnote{Ver, por exemplo, \citeonline{Parsons1984} e \citeonline{Burge1984}.}.
O paradoxo de Kripke, entretanto, deixa claro que não é este o único
problema com as hierarquias de Tarski, mas também o fato de seus níveis
serem determinados por aspectos intrínsecos das sentenças.

O ponto aqui é que tais hierarquias seriam absolutamente insensíveis
às circunstâncias empíricas por meio das quais determinadas sentenças
podem engendrar paradoxos. Note que não é possível estabelecer sintaticamente
os níveis de (15) e (16), na medida em que eles se sobrepõem mutuamente.
O predicado verdade da sentença (15) deveria pertencer a um nível
superior ao da sentença (16) e o predicado da sentença (16) deveria
pertencer a um nível superior àquele da sentença (15).

Além disso, é certamente possível que as circunstâncias favorecessem
atribuições consistentes de valores a (15) e (16). Suponha, por exemplo,
que Jones asserisse mais alguma sentença reconhecidamente falsa sobre
o Watergate. Neste caso, a maioria das asserções feitas por Nixon
a respeito são, de fato, falsas. Podemos, assim, admitir consistentemente
que (15) é verdadeira e (16) é falsa. Um alteração feita sobre as
circunstâncias, aspectos extrínsecos às sentenças, dissolve o paradoxo.
Isso indicaria que os paradoxos da linguagem natural podem depender
de aspectos extrínsecos às sentenças.

As hierarquias de inspiração tarskiana continuam tendo um papel importante
na literatura atual em torno dos paradoxos, mesmo que a objeção de
Kripke tenha aberto um novo ramo de trabalhos e abordagens. O ponto
essencial que podemos extrair aqui, portanto, é que, a despeito do
diagnóstico clássico se apoiar em resultados precisos e aprofundar
a compreensão dos problemas envolvidos nos paradoxos, não é claro
qual seja a lição que podemos extrair dele a respeito do conceito
ordinário de verdade e do significado do paradoxo do mentiroso na
linguagem natural. Essas ainda são questões em aberto.

\section{Considerações finais}

A posição deflacionista pode ser compreendida não como a afirmação
de que nada há para ser dito acerca da verdade além do que é dito
por instâncias do esquema T {em geral}, mas sim no {âmbito de uma
investigação lógico-filosófica.} Sentenças T poderiam, por assim
dizer, ser aperfeiçoadas, mas isso não seria um problema filosófico,
nem lógico. Fornecer uma resposta {substantiva} acerca das condições
de verdade de uma sentença não seria tarefa da filosofia, mas sim
do setor específico do conhecimento que trata a sentença em questão.
Assim, a tarefa de complementar o que é afirmado pelas instâncias
de T não pertence à filosofia. Quem pode fornecer informações não
triviais acerca da verdade da sentença `NaCl dissolve na água' não
é o filósofo, mas sim o químico, do mesmo modo que é o matemático
que vai dizer algo não trivial acerca da falsidade da sentença `51
é um número primo'. Não era intenção de Tarski entrar em detalhes
acerca das condições de verdade de uma dada sentença, mas apenas {expressar}
tais condições de maneira correta.\footnote{``De fato, a definição semântica  da verdade não implica nada a respeito de condições sob as quais uma sentença como (1):\\ (1) \textit{a neve é branca}\\ possa ser afirmada. Ela implica apenas que, em quaisquer circunstâncias em que afirmemos ou neguemos essa sentença, devemos estar prontos para afirmar ou negar a sentença correlata (2)\\ (2) \textit{a sentença `a neve é branca' é verdadeira}'' (CSV, p. 189).} E de fato, instâncias de T expressam corretamente as condições de
verdade de uma dada sentença, independentemente de fazê-lo utilizando
a própria sentença. Nessa perspectiva, o ponto não é considerá-las
triviais, pois trata-se de uma posição filosófica acerca da verdade:
dizer mais sobre a verdade de uma sentença não seria tarefa da filosofia.

Portanto, é bastante razoável afirmar que a teoria de Tarski diz tudo
o que Tarski achava que havia para ser dito, no âmbito de uma investigação
lógico-filosófica, acerca da relação de correspondência entre uma
sentença e a realidade. E isso não deixa de ser uma posição filosoficamente
relevante sobre o problema da verdade.

De outro lado, os resultados formais (positivos e negativos) acerca
do conceito de verdade determinaram importantes desdobramentos em
Lógica. O método utilizado por Tarski na definição recursiva de satisfação
é o ponto de partida para as semânticas formais atualmente utilizadas
em Lógica de Predicados. Ademais, o Teorema da Indefinibilidade de
Tarski deu origem a um novo campo de pesquisas, aquele das teorias
axiomáticas da verdade\footnote{Um importante \textit{handbook} em teorias axiomáticas, \citeonline{Horsten2011},
leva o título \textit{The Tarskian Turn}.}.

\nocite{Chateaubriand2001}

\nocite{Fenstad2004}

\nocite{Frost-Arnold2004}

\nocite{Heck2007}

%
%
%

\nocite{Sher1999a}

\nocite{Sher2001}

\nocite{Tarski1969}

%

\bibliographystyle{abntex2-alf}
\bibliography{abilio}

\end{document}